\documentclass[11pt,leqno]{article}
\usepackage{graphicx, amsfonts, amsthm, amsxtra, verbatim, multicol}
\usepackage[mathscr]{euscript}

\begin{document}
\title{
The infinitesimal 16th Hilbert problem in dimension zero
}
\author{Lubomir Gavrilov \\
\normalsize \it   Universit\'e
Paul Sabatier\\
\normalsize \it Laboratoire Emile Picard, MIG, CNRS UMR 5580\\
\normalsize \it  31062 Toulouse Cedex 9, France \\
\\ Hossein Movasati\\
\normalsize \it Technische Universit\"{a}t Darmstadt\\
\normalsize \it Fachbereich Mathematik\\
\normalsize \it Schlossgartenstr. 7, 64289 Darmstadt, Germany } \date{July 1, 2005}
\maketitle
\begin{abstract}
We study the analogue of the infinitesimal 16th Hilbert problem in dimension zero. Lower and upper bounds for
the number of the zeros of  the corresponding Abelian  integrals  (which are algebraic functions) are found. We
study the relation between the vanishing of an Abelian integral $I(t)$ defined over $\mathbb{Q}$ and  its
arithmetic properties. Finally, we give necessary and sufficient conditions for an Abelian integral  to be
identically zero.

\end{abstract}

\def\Z{\mathbb{Z}}                   
\def\F{{\cal   F}}                   
\def\A{\mathbb{A}}                   
\def\G{\mathbb{G}}                   
\def\K{\mathsf{k}}                   
\def\Q{\mathbb{Q}}                   
\def\C{\mathbb{C}}                   
\def\R{\mathbb{R}}                   
\def\N{\mathbb{N}}                   
\def\Sc{\mathsf{S}}                  
\def\P{\mathbb{P}}
\def\Gal{{\mathrm Gal}}              
\newcommand{\mat}[4]{
     \begin{pmatrix}
            #1 & #2 \\
            #3 & #4
       \end{pmatrix}
    }                                
\newtheorem{theo}{Theorem}
\newtheorem{exam}{Example}
\newtheorem{coro}{Corollary}
\newtheorem{defi}{Definition}
\newtheorem{prob}{Problem}
\newtheorem{lemm}{Lemma}
\newtheorem{prop}{Proposition}
\newtheorem{rem}{Remark}
\newtheorem{conj}{Conjecture}
\section{Introduction}
Let $f: X \rightarrow Y$ be a morphism of complex algebraic varieties which defines a locally trivial
topological fibration. Let $\gamma(a) \in H_q(f^{-1}(a),\Z)$ be a continuous family of $q$-cycles and $\omega $
be a regular $q$-form on $X$ which is closed on each fiber $f^{-1}(a)$ (the latter is always true if $q= \dim
f^{-1}(a)$). By Abelian integral (depending on a parameter) we mean a complex multivalued function of the form
$$
I: Y \rightarrow \C : a \rightarrow I(a)= \int_{\gamma(a)} \omega .
$$
Through the paper we shall also suppose that the varieties $X,Y$, the morphism $f$ and the $q$-form $\omega$ are
defined over a subfield $\K \subset \C$. In the case
$$\K = \R,\  Y = \C \setminus S,\  X=f^{-1}(Y)\subset \C^2 , \
\omega =P\,dx + Q\,dy, \ f,P,Q \in \R[x,y]
$$
where $S$ is the finite set of atypical values of $f$,  the zeros of $I(a)$ on a suitable open real interval are
closely related to  the limit cycles of the perturbed foliation on the real plane $\R^2$ defined by
$$
df + \epsilon \omega=0,\  \epsilon\sim 0 .
$$
Recall that the second part of the 16th Hilbert problem asks to determine the maximal number and positions of
the limit cycles of a polynomial plane vector field (or foliation) of a given degree. The infinitesimal 16th
Hilbert problem asks then to \emph{find the
 the exact
upper bound $Z(m,n)$ for the number of the zeros of $I(a)$ on an open interval, where $\deg f \leq m$, $\deg P,
Q \leq n$} \cite{ily1}. It is only known that $Z(m,n) < \infty$ \cite{va,kh} and $Z(3,2)=2$ \cite{Gav4}.

More generally, let $X$ and $Y$ be Zariski open subsets in $\C^{q+1}$ and $\C$ respectively, $f$ a polynomial
and $\omega$ a polynomial $q$-form in $\C^{q+1}$,
all these objects being defined over a subfield $\K \subset
\C$. \emph{What is the exact upper bound $Z(m,n,\K,q)$ for the number of
the zeros $a\in \K \cap {\cal D}$ of the
Abelian integral $I$?} Here ${\cal D}$ is any simply connected domain in $Y$.

The present paper addresses the above question in the simplest case $q=0$. The Abelian integral $I$ is then an
algebraic function over $\K[a]$ and every algebraic function over $\K[a]$ is an Abelian integral defined over
$\K$. We prove in Theorem \ref{main} that

\begin{equation}\label{ineq}
n-1- [\frac{n}{m}] \leq Z(m,n,\K,0) \leq \frac{(m-1)(n-1)}{2} .
\end{equation}
The lower bound in this inequality is given by the dimension of the vector space of Abelian integrals
$$ V_n= \{ \int_{\gamma(a)} \omega, \  \deg \omega \leq n\}$$
where $f$ is a fixed general polynomial of degree $m$, while the upper bound is a reformulation of Bezout's
theorem.  When $d=3$ we get $ Z(d,d-1,\K,0)=1$. We give some evidence in Proposition \ref{p6}  that, in the case
$\K=\R$, $m=4$, $n=3$, the upper bound of (\ref{ineq}) is strictly bigger than $ Z(4,3,\R,0)$. This Proposition
also suggests that
$$
\lim_{d\rightarrow \infty} \frac{Z(d,d-1,\R,0)}{d} = 1
$$
or, in other words, the space of Abelian integrals $V_d$ is Chebishev, possibly with some accuracy. Recall that
$V_n$ is said to be Chebishev with accuracy $c$ if every $I\in V_n$ has at most $\dim V_n -1 +c$ zeros in the
domain $D$. The Chebishev property (if satisfied) would mean that the infinitesimal 16th Hilbert problem in
dimension zero is a problem of real algebraic geometry (as opposed to Bezout's theorem which is a result of
complex algebraic geometry).

To the rest of the paper we explore some arithmetic properties of  Abelian integrals.
 When $\K$
is a number field and $q=1$ (so $f^{-1}(a)$ is a smooth curve) a recent theorem of Movasati and Viada
\cite{evho} gives necessary and sufficient conditions for the vanishing of an Abelian integral $I(a)$ defined
over $\K$. We formulate  the 0-dimensional analogue of this result (Theorem \ref{3may05}). Its proof uses the
relation between the vanishing of an Abelian integral of dimension zero and the Galois group of the splitting
field of $f-a$. Finally we make use of the monodromy of $f$ to obtain two additional results. The first one
improves the upper bound for the number of the zeros of an Abelian integral with fixed $f$ (Theorem
\ref{main1}). The second one gives necessary and sufficient conditions for an Abelian integral $I(a)$ to be
\textit{identically} zero. The analogue of this result for $q=1$ is not known, although it is essential for
computing the so called higher order Poincar\'{e}-Pontryagin functions \cite{gi,gav5}.

The paper is organized as follows. In sections \ref{s2} and \ref{s3} we summarize, for convenience of the
reader, the basic properties of zero dimensional Abelian integrals and their algebraic counterpart: the global
Brieskorn module. The canonical connection of the (co)homology bundle is explained on a simple example of a
 polynomial $f$ of degree three. In section \ref{s4} we prove the Bezout type estimate for
$Z(m,n,\K,0)$ and consider the examples $\K=\R$, $m=3,4$. The arithmetic aspects of the problem are treated in
section \ref{s5}, and the monodromy group of $f$ in section \ref{s6}.

{\bf Acknowledgments.} The main results of the paper were obtained
 while the second author was visiting the University Paul Sabatier of
Toulouse. He thanks for its hospitality.

\section{Zero dimensional Abelian integrals }
\label{s2}

In this section we introduce the necessary notations and prove, for convenience of the reader, some basic facts
about the Abelian integrals of dimension zero.

Let $M= \{x_1,x_2,\dots,x_d\}$ be a discrete topological space and $G$ an additive abelian group. By abuse of
notation we denote by $H_0(M,G)$ the reduced homology group
$$
H_0(M,G) = \{ \sum_{i=1}^d n_i x_i : n_i \in G, \sum_{i=1}^d n_i = 0\} .
$$
It is a free $G$-module of rank $d-1$ generated by
$$
x_1-x_d, x_2-x_d, \dots , x_{d-1}-x_d
$$
and its dual space will be denoted by $H^0(M,G)$.
To the polynomial
$$
f(x,a)= x^d-a_1 x^{d-1} -\cdots -a_d, a=(a_1 , a_2, \dots ,a_d)
$$
we associate the surface
$$
V= \{(x,a)\in \C^{d+1}: f(x,a) = 0 \}
$$
and the (singular) fibration
\begin{equation}
\label{fibration}
V \rightarrow \C^d : (x,a)\mapsto a
\end{equation}
with fibers $L_a= \{x\in \C: f(x,a)=0\}$. The polynomial $f(x,a)$ is a versal deformation of the singularity
$f(x,0)=x^d$ of type $A_{d-1}$.
 We denote by $\triangle(a)$ the discriminant of $f(x,a)$ with respect
to $x$. The corresponding discriminant locus $\Sigma = \{a\in \C^d: \triangle(a) = 0 \}$ is the set of
parameters $a$, such that $f(x,a)$ has a multiple root (as a polynomial in $x$).

The map (\ref{fibration}) induces homology and co-homology bundles with base $\C^d \setminus \Sigma$, and fibers
$H_0(L_a,\Z)$ and $H^0(L_a,\C)$ . The continuous families of cycles
$$\gamma_{ij}(a)=x_i(a)-x_j(a) \in H_0(L_a,\Z)$$
 generate a basis of locally constant
sections of a unique connection in the homology bundle (the so called Gauss-Manin connection).

Let $\K\subset \C$ be a field. To define the connection algebraically we need the global Brieskorn module
(relative co-homology) of $f$ which is defined as
\begin{equation}\label{brieskorn}
H = \frac{\K[x,a]}{ f. \K[x,a] + \K[a]} .
\end{equation}
This is  a $\K[a]$-module in an obvious way. The basic properties of such modules in the local multi-dimensional
case ($x\in \C^{n+1}, a\in \C$) when $\K[a]$ is replaced by $\C\{a\}$ were studied by Brieskorn \cite{bri} and
Sebastiani \cite{seb}. The first results in the global one-dimensional case ($x\in \C^2$) were proved in
\cite{gav}. For  arbitrary $n$ see \cite{dim}, \cite{sab}, \cite{bon} and \cite{mo}. In the zero-dimensional
case the main properties of $H$ are rather obvious and are summarized in Proposition \ref{prop1} and Proposition
\ref{prop2} bellow.

\begin{prop}
\label{prop1} $H$ is a free $\K[a]$-module of rank $d-1$ generated by
$x, x^2, \dots , x^{d-1}$. More precisely,
for every $m\geq d$ the following identity holds in $H$
$$
x^m = \sum_{i=1}^{d-1} p_i(a) x^i
$$
 where $p_i(a)\in \K[a]$ are suitable  weighted homogeneous
polynomials of degree  $m-i$, and $weight(a_i)=i$.
\end{prop}
\begin{proof}
The proof is by induction on $m$ and is  left to the reader.
\end{proof}
Let $\gamma(.)$ be a locally constant section of the homology bundle of $f$. Every $\omega= \omega(x,a) \in
\K[x,a]$ defines a global section of the co-homology bundle by the formula
\begin{equation}\label{integral}
I(a) = \int_{\gamma(a)} \omega = \sum_{i}\omega (x_i(a),a)
\end{equation}
where $\gamma(a)= \sum_{i}n_i x_i(a)$, $\sum_{i} n_i=0$.
\begin{defi}
 An Abelian integral of dimension zero over the field $\K$ is a function $I(a)$, $a\in \C^d$, of the form
(\ref{integral}), where $f, \omega \in \K[x,a]$ and $\gamma(a)$ is a continuous family of cycles.
\end{defi}
\begin{rem}
A ( multivalued ) function $I(a)$  is an Abelian integral if and only if it is an algebraic function. Indeed,
let $x=x(a)$ be the algebraic function defined by $g(x,a) \equiv 0$, $g\in  \K[x,a]$. Then it is an Abelian
integral $I(a)$ defined by either
$$
\omega = x, f= x g(x,a), \gamma(a) = x(a)- 0
$$
or
$$
\omega = x, f=g(2x,a)g(-2x,a), \gamma(a) = \frac{x(a)}{2}-
(-\frac{x(a)}{2}) .
$$
In the next two propositions we shall suppose, however, that
$$
f(x,a)= x^d-a_1 x^{d-1} -\cdots -a_d, a=(a_1 , a_2, \dots ,a_d) .
$$
\end{rem}
\begin{prop}
\label{prop2} The polynomial $\omega= \omega(x,a) \in \K[x,a]$ defines the zero section of the canonical
co-homology bundle of $f$, if and only if $\omega$ represents the zero equivalence class  in the global
Brieskorn module (\ref{brieskorn}).
\end{prop}
\begin{proof}
Indeed, if $[\omega]= 0$ in $H$,
the claim is obvious. If $\omega$ defines the zero section of the co-homology
bundle, then
$$
\omega (x_i(a),a) = \omega (x_j(a),a), \forall a,i,j .
$$
According to Proposition \ref{prop1} we may suppose that
$$
\omega (x,a)= \sum_{i=1}^{d-1} p_i(a) x^i + f.p(x,a) + q(a)
$$
and  $a$ is such that $x_i(a)\neq x_j(a)$ for $i\neq j$. Then it follows that the affine algebraic curves
\begin{equation}\label{c1}
\Gamma_\omega = \{(x,y)\in \C^2: \frac{\omega (x,a) - \omega (y,a)}{x-y}=0 \}
\end{equation} and
\begin{equation}\label{c2}
\Gamma_f = \{ (x,y)\in \C^2: \frac{f (x,a) - f(y,a)}{x-y} =0\}
\end{equation}
have at least $d(d-1)$ distinct intersection points at the points$(x,y)=(x_i,x_j)$. But this contradicts the
Bezout's theorem, as the degree of the curves (\ref{c1}) and (\ref{c2}) is (at most) $d-2$ and $d-1$
respectively. It follows that either $\Gamma_f$ and $\Gamma_\omega$ have a common component, or $p_i(a)=0$,
$\forall a,i$. In the former case the algebraic curve is reducible which is impossible for generic values of
$a$. We obtain finally that $\omega(x,a)$ represents the zero equivalence class in the Brieskorn module
$H$.
\end{proof}
For a given $f\in \K[x,a]$ and a section $\gamma(a)\in H_0(L_a,\Z)$,
let ${\cal A}_f$ be the set of
Abelian integrals (\ref{integral}), where $a$ belongs to some simply
 connected sub-domain of $\C$. Then ${\cal
A}_f$  is a $\K[a]$-module and moreover
\begin{prop}
\label{prop3} ${\cal A}_f$  and the Brieskorn module $H$ are  isomorphic $\K[a]$-modules.
\end{prop}
\begin{proof}
The homomorphism
$$
H \rightarrow {\cal A}_f : \omega \rightarrow \int_{\gamma(a)} \omega = \sum_{i}\omega (x_i(a),a)
$$
is obviously surjective. As the monodromy group of the fibration defined by $f$ is transitive, then
 $I(a)\equiv 0$ implies that
  $\int_{\gamma(a)} \omega \equiv 0$ for every section $\gamma(a) \in H_0(f^{-1}(a),\Z)$. Proposition
  \ref{prop2} implies that $\omega=0 \in H$.
  \end{proof}

Let $\K^q \rightarrow \K^d: b\rightarrow a$ be a polynomial map,
and consider
\begin{equation}\label{brieskorn1}
H_b = \frac{\K[x,b]}{ g. \K[x,b] + \K[b]}
\end{equation}
where
$$
g= x^d -a_1(b) x^{d-1}- \dots - a_d(b) .
$$
As before $H_b$ is a free $\K[b]$ module with generators
$x, x^2,\dots,x^{d-1}$.
Of particular interest is the polynomial map
$$
\K \rightarrow \K^d : t \rightarrow (a_1^0,a_2^0,\dots,a_{d-1}^0,t) .
$$
The $\K[t]$ module $H_t$ is then isomorphic to
\begin{equation}\label{brieskorn2}
H_t = \frac{\K[x]}{ \K[g]}, g=  x^d -a_1^0 x^{d-1}- \dots - a_{d-1}^0x
\end{equation}
with multiplication $t\cdot \omega = g(x)\omega(x) \in \K[x]$. The analogues of Proposition \ref{prop1},
\ref{prop2} hold true for $H_t$, but not Proposition \ref{prop3} (see section \ref{monodromy}).

\section{The connection of $H$}
\label{s3}
 \label{connection} Let $x(a)$ be a root of the polynomial $f$, $f(x(a),a)\equiv 0$. Then
\begin{equation}\label{eq1}
\frac{\partial x(a)}{\partial a_i} \frac{\partial f(x,a)}{\partial x} \equiv x(a)^i .
\end{equation}
There exist polynomials $p,q\in \Q[x,a]$ such that
\begin{equation}\label{eq2}
p \frac{\partial f}{\partial x} + q f = \triangle(a)
\end{equation}
and hence in $H$
$$\triangle(a)  = p(x,a) \frac{\partial f}{\partial x} .$$
This combined with (\ref{eq1}) suggests to define a connection on $H$ as follows
$$
\nabla_{\frac{\partial }{\partial a_i}}
 : H \rightarrow H_{\Delta}: x^m  \mapsto \frac{m x^{m-1+i}p(x)}{\Delta}
$$
where  $H_\Delta$ is the localization of $H$ on $\{\Delta,\Delta^2,\cdots\}$. The operator $\nabla$ satisfies
the Leibniz rule and so it is a connection on the module $H$. It follows from (\ref{eq1}), (\ref{eq2}) that
\begin{equation}\label{eq3}
\frac{\partial }{\partial a_i}\int_{\gamma(a)}\omega=\int_{\gamma(a)}\nabla_{\frac{\partial }{\partial
a_i}}\omega
\end{equation}
where $\gamma(a)$ is a continuous family of  cycles. Every element of $H$ defines a section of the co-homology
bundle of $f$. By (\ref{eq3}) every continuous family of  cycles is a locally constant section of the homology
bundle, which means that $\nabla$ coincides with the Gauss-Manin connection described previously

It is well known that $\nabla$ is a flat (integrable) connection. Indeed, a fundamental matrix of solutions for
this connection is given by
$$
X(a) = ( \int_{\gamma_j(a)} x^i)_{i,j=1,1}^{d-1,d-1}
$$
where $\gamma_1(a), \gamma_2(a),\dots , \gamma_{d-1}(a)$ is a basis of locally constant sections of the homology
bundle. The connection form is therefore
$$
 \sum_{i=1}^d A_i(a) d a_i , \mbox{  where  } A_i(a)=\frac{\partial }{\partial a_i} X(a) . X^{-1}(a)
$$
and the family of (Picard-Fuchs) differential operators
$$
\frac{\partial }{\partial a_i}- A_i , i=1,2,\dots ,d
$$
commute. We end this section by a simple but significant example. Let
$$
f= 4 x^3 - g_2 x - g_3
$$
with discriminant
$$
\triangle(g_2,g_3) = g_2^3- 27 g_3^2 .
$$
Suppose further that $g_2, g_3$ depend on a parameter $z$. A straightforward and elementary computation implies
\begin{prop}
In the Brieskorn module $H$ the following identity holds
$$
\triangle(z) \nabla_{\frac{\partial}{dz} }
\left(%
\begin{array}{c}
  x \\
  x^2 \\
\end{array}%
\right)
=
\left(%
\begin{array}{cc}
  \frac{\triangle'_z}{6} & -3 \delta \\
  - \frac{g_2 \delta}{2} & \frac{\triangle'_z}{3}\\
\end{array}%
\right)
\left(%
\begin{array}{c}
  x \\
  x^2 \\
\end{array}%
\right)
$$
where
$$
\delta(z) = 3 g_3 \frac{d g_2}{dz} - 2 g_2 \frac{d g_3}{dz} .
$$
If we introduce the Abelian integrals
$$
\left(%
\begin{array}{c}
  \eta_1 \\
  \eta_2 \\
\end{array}%
\right)
=
\left(%
\begin{array}{c}
  \frac{\int_{\gamma(z)} x}{\Delta^{1/4}} \\
  \frac{ \int_{\gamma(z)} x^2 }{\Delta^{1/4} }\\
\end{array}%
\right)
$$
then they satisfy the Picard-Fuchs system
\begin{equation}\label{pf}
\triangle(z) \frac{d}{dz}
\left(%
\begin{array}{c}
  \eta_1 \\
  \eta_2 \\
\end{array}%
\right) =
\left(%
\begin{array}{cc}
  -\frac{\triangle'_z}{12} & -3 \delta \\
  - \frac{g_2 \delta}{2} & \frac{\triangle'_z}{12}\\
\end{array}%
\right)
\left(%
\begin{array}{c}
  \eta_1 \\
  \eta_2 \\
\end{array}%
\right)
\end{equation}
\end{prop}

It is interesting to compare the above  to the Picard-Fuchs system associated to the ''stabilization'' $y^2- 4
x^3 + g_2(z) x + g_3(z)$ of $f$. Namely, let
$$
\eta_1= \int_{\gamma(z)} \frac{dx}{y}, \ \eta_2= \int_{\gamma(z)} \frac{xdx}{y}
$$
be complete elliptic integrals of first and second kind on the elliptic curve with affine equation
$$
\Gamma_z = \{(x,y)\in \C^2: y^2= 4 x^3 - g_2(z) x - g_3(z) \}
$$
where $\gamma(z) \subset \Gamma_z$ is a continuous family of closed loops (representing a locally constant
section $z\mapsto H_1(\Gamma_z,\Z)$ of the homology bundle). Then $\eta_1, \eta_2$ satisfy the following
Picard-Fuchs system
(this goes back at least to \cite[Griffiths]{gri}, see \cite[Sasai]{sas})
\begin{prop}
\begin{equation}\label{pf2}
\triangle(z) \frac{d}{dz}
\left(%
\begin{array}{c}
  \eta_1 \\
  \eta_2 \\
\end{array}%
\right) =
\left(%
\begin{array}{cc}
  -\frac{\triangle'_z}{12} & -\frac{3 \delta}{2} \\
  - \frac{g_2 \delta}{8} & \frac{\triangle'_z}{12}\\
\end{array}%
\right)
\left(%
\begin{array}{c}
  \eta_1 \\
  \eta_2 \\
\end{array}%
\right) .
\end{equation}
\end{prop}
Note the striking similarity of these two non-equivalent systems.
The algorithms which calculate the Gauss-Manin connection can be implemented
in any software for commutative algebra
(see \cite{mo05}). Using them one can obtain equalities like (\ref{pf2}) for
other families of varieties in arbitrary dimension.
\section{The infinitesimal 16th Hilbert problem in dimension zero}
\label{s4}
 It is well known that the  the number of the limit cycles of the perturbed real foliation
$$
df + \varepsilon (P dx + Q dy) = 0,\  P,Q,R \in \R[x,y],\  \varepsilon \sim 0
$$
is closely related to the number of the zeros of the Abelian integral (Poincar\'{e}-Pontryagin function)
$$
I(t)= \int_{\gamma(t)} P dx + Q dy
$$
where $\gamma(t)\in H_1(f^{-1}(t),\Z)$ is a continuous family of cycles. The problem on zeros of such Abelian
integrals, in terms of the degrees of $f,P,Q$, is known as the infinitesimal 16th Hilbert problem: see the
recent survey of Ilyashenko \cite{ily1},
as well \cite[problem 1978-6]{ar05}. This problem is still open (except
in the case $\deg F \leq 3$, see \cite{Gav4,hi98}). One can further generalize, by taking $f\in
\R[x_1,x_2,\dots,x_n]$, $\omega $ - a polynomial $n-1$ form, $\gamma(t) \in H_{n-1}(f^{-1}(t),\Z)$ a locally
constant section of the homology bundle of $f$, and
$$
I(t) = \int_{\gamma(t)} \omega .
$$
In the present paper we solve (partially) the infinitesimal  16th Hilbert problem by taking $n=1$. The Abelian
integral $I(t)$ is of dimension zero, in the sense explained in the preceding section. To our knowledge, such
integrals appeared for a first time, in  the context of the 16th Hilbert problem, in the Ilyashenko's pioneering
paper \cite{ily}, see \cite{gav}.

To formulate the problem, let us denote $f \in \K[x]$ where $\K\subset \C$ is a field, and consider the singular
fibration
$$
f: \C \rightarrow \C : x \mapsto f(x)
$$
with fibers $L_t = f^{-1}(t)$. Let ${\cal D} \subset \C\setminus \Sigma$ be simply connected set, where $\Sigma$
is the set of critical values of $f$. A cycle $\gamma(t) \in H_0(L_t,\Z)$ is said to be simple, if $\gamma(t)=
x_i(t)-x_j(t)$, where $f(x_i(t))= f(x_j(t)) = t$. Let $\gamma(t)$ be a continuous family of simple cycles, $m,n$
two integers and $\K$ a field.
The Infinitesimal 16th Hilbert problem in dimension zero is

 \textit{Find the exact upper bound $Z(m,n,\K,0)$ of the number of
the zeros
$$
\{ t\in \K\cap {\cal D}: I(t)=0,\  \deg f \leq m,\  \deg \omega \leq n,
\ {\cal D} \subset \C\setminus \Sigma  \}
$$
 where  ${\cal D}$ is any simply connected complex domain. }

 In the present paper we are interested in the cases
$\K=\Q,\R, \C$. When $\K=\R$ the problem has the
following geometric interpretation. Consider the real plane algebraic curves
$$
\Gamma_f^\R = \{(x,y)\in \R^2: \frac{f(x)-f(y)}{x-y}= 0\},\  \Gamma_\omega^\R = \{(x,y)\in \C^2:
\frac{\omega(x)-\omega(y)}{x-y}= 0\}.
$$
The Abelian integral
$$
\int_{\gamma(t)} \omega = \omega (x_i(t))-\omega (x_j(t))
$$
vanishes if and only if $(x_i(t),x_j(t))\in \Gamma_f^\R \cap \Gamma_\omega^\R$. If we suppose in addition that
${\cal D} \subset \R$ is an open interval, then $(x_i(t),x_j(t)): t\in {\cal D}$ is contained in some connected
component of $\Gamma_f^\R$ which we  denote  by $\Gamma_{f,0}^\R$.

It is clear that the number of intersection points $\#(\Gamma_{f,0}^\R,\Gamma_\omega^\R)$ (counted with
multiplicity)
 between $\Gamma_{f,0}^\R$ and the real algebraic curve
$\Gamma_\omega^\R$ is an upper  bound for the corresponding number $Z(m,n,\R,0)$. On the other hand
$\#(\Gamma_{f,0}^\R,\Gamma_\omega^\R)$ can be bounded by the Bezout's theorem. It is not proved, however, that
$$
Z(m,n,\R,0) = \#(\Gamma_{f,0}^\R,\Gamma_\omega^\R)
$$
and we discuss this at the end of the section.  We have the following
\begin{theo}
\label{main}
\begin{equation}
\label{maineq}
n-1- [\frac{n}{m}] \leq Z(m,n,\C,0) \leq \frac{(m-1)(n-1)}{2} .
\end{equation}
\end{theo}
\begin{proof}
Let $\gamma(t) = x(t) -y(t)$ be a continuous family of simple cycles. Then $I(t)=0$ for some
$t\not\in \Sigma$ if and only if $(x(t),y(t))$ is an isolated intersection point of the plane algebraic curves
$$
\Gamma_f = \{(x,y)\in \C^2: \frac{f(x)-f(y)}{x-y}= 0\}, \Gamma_\omega = \{(x,y)\in \C^2:
\frac{\omega(x)-\omega(y)}{x-y}= 0\}.
$$
Indeed, if an intersection point were non isolated, this would mean that the curves $\Gamma_f$ and
$\Gamma_\omega$ have a common connected component and $I(t) \equiv 0$.
 By Bezout's Theorem the number of isolated intersection points of $\Gamma_f$ and
$\Gamma_\omega$, counted with multiplicity, is bounded by $(\deg f - 1)(\deg \omega -1)$. Moreover if $(x,y)$ is
an isolated intersection point which corresponds to some $t_0 \in \C\setminus \Sigma$, $\gamma(t_0) = x-y$,
$I(t_0)=0$, then $(y,x)$ is an isolated intersection point too. As $I(t)$ is single-valued in ${\cal D}$ then
$(y,x)$ does not correspond to any zero of $I(t)$ in ${\cal D}$. Thus the number of the zeros of $I(t) $ on
${\cal D}$ is bounded by $(\deg f - 1)(\deg \omega -1)/2$.

Let
$$ V_n= \{ \int_{\gamma(t)} \omega, \ \deg \omega \leq n\}$$
be a vector space of Abelian integrals defined in a simply connected domain for some fixed \emph{generic}
polynomial $f$ of degree $m$. We have $\dim V_n -1 \leq Z(m,n,\C,0)$. On the other hand, if $f$ is a generic
polynomial, then the orbit of $\gamma(t)$ under the action of the monodromy group of the polynomial spans
$H_1(f^{-1}(t),\Z)$ (see sectiuon \ref{s6}). Therefore $I\in V_n$ is identically zero if and only if $\omega$
represents the zero co-homology class in $H^1(f^{-1}(t),\C)$ which is equivalent (by Proposition \ref{prop2}) to
$\omega \in \C[f]$. This shows that the vector space $V_n$ is isomorphic to
$$\{ \omega\in \C[x]: \deg \omega \leq n\} / \{ p(f(x)): p \in \C[x]: \deg p(f(x)) \leq n\} .$$
The basis of this space is
$$
\{x^i f^j: i+j\,m \leq n \}
$$
and hence $\dim V_n = n- [\frac{n}{m}]$. The Theorem is proved.
\end{proof}

The bound in the above Theorem is probably far from the exact one.
If one wants to count zeros in a not simply connected domain ${\cal D}$
then the bound is exact for the case $\deg(\omega)=d-1$. For
instance take $f=\Pi_{i=1}^d(x-i)$ and $\omega=\Pi_{i=1}^{d-1}(x-i)$. Then
$\int_{\gamma_{ij}}\omega=0$ for all
$\gamma_{ij}=i-j\in H_0(\{f=0\},\Z), i,j=1,2,\ldots,d-1$.
In \S \ref{monodromy} we will give another approach
using the monodromy  representation of $f$.
\begin{exam}
Let $\deg f = 3$, and  $V$ be the $\K$-vector space of Abelian integrals generated by
$$
\int_{\gamma(t)} x, \int_{\gamma(t)} x^2 .
$$
By Theorem \ref{main} each $I \in V$ has at most one simple zero in ${\cal D}$. As $\dim V = 2$ this means that
the bound can not be improved, that is to say  $V$ is a Chebishev space in ${\cal D}$.
\end{exam}
\begin{exam}
Let $\deg f = 4$, and  $V$ be the $\K$-vector space of Abelian integrals generated by
$$
\int_{\gamma(t)} x, \int_{\gamma(t)} x^2 , \int_{\gamma(t)} x^3.
$$
By Theorem \ref{main} each $I \in V$ has at most three zeros
in ${\cal D}$. As $\dim V = 3$ this does not imply
that $V$ is Chebishev.
\end{exam}

\begin{exam}

Consider the particular case $f= x^4-x^2$ and take $\K= \R$. The set of critical values is $\Sigma = \{0,-1/4\}$
and let  $\gamma(t)$ be a continuous family of simple cycles  where $t\in (-1/4,0)$.
\begin{prop}
\label{p6}
 The Vector space $V_n$ of Abelian integrals $I(t)=\int_{\gamma(t)}\omega$, $\deg \omega \leq n$ is
Chebishev. In other words, each $I\in V_n$ can have at most $\dim V_n -1$ zeros in $(-1/4,0)$.
\end{prop}
\begin{proof}
We shall give two distinct proofs. The first one applies only for the family of cycle vanishing at $t=0$, but
has the advantage to hold in a complex domain.

1) Suppose that $ \gamma(t)= x_1(t)-x_2(t)$ is a continuous family of simply cycles vanishing as $t$ tends to
zero and defined for $t \in \C\setminus (-\infty, -1/4]$. Each integral  $I\in V_n$ admits  analytic
continuation in $\C\setminus (-\infty, -1/4]$. Following Petrov \cite{Petrov86}, we shall count the zeros of
$I(t)$ in $\C\setminus (-\infty, -1/4]$ by making use of the argument principle. Consider the function
$$
F(t)=   \frac{\int_{\gamma(t)} \omega}{\int_{\gamma(t)} x}
$$
which admits analytic continuation in $\C\setminus (-\infty, -1/4]$. Let $D$ be the domain obtained from
$\C\setminus (-\infty, -1/4]$ by removing a "small" disc $\{z\in\C: |z|\leq r\}$ and a "big" disc $\{z\in\C:
|z|\geq R \}$. We compute the increase of the argument of $F$ along the boundary of $D$ traversed in a positive
direction. Along the boundary of the small disc the increase of the argument is close to zero or negative
(provided that $r$ is sufficiently small). Along the boundary of the big disc the increase of the argument is
close to $(n-1)\pi/2$ or less than $(n-1)\pi/2$. Finally, along the coupure  $(-\infty, -1/4)$ we compute the
imaginary part of $F(t)$. Let $F^\pm(t),\gamma^\pm(t)$ be the determinations of $F(t), \gamma(t)$ when
approaching $t\in (-\infty, -1/4)$ with $Im t
>0$ ($Im t< 0$). We have
$$
Im F(t) = (F^+(t)-F^-(t))/2\sqrt{-1}, \gamma^+(t)-\gamma^-(t)=\delta(t)
$$
where $\delta(t)= x_3(t)-x_1(t)-(x_4(t)-x_2(t))$ and $x_3(t)$, $x_4(t)$ are roots of $f(x)-t$ which tend to
$x_1(t)$ and $x_2(t)$ respectively, as $t$ tends to $-1/4$ along a path contained in $\C\setminus (-\infty,
-1/4]$. We obtain
$$
2\sqrt{-1} Im F(t) = \beta \frac{\det
\left(%
\begin{array}{cc}
  \int_{\gamma(t)} \omega & \int_{\gamma(t)} x \\
 \int_{\delta(t)} \omega & \int_{\delta(t)} x \\
\end{array}%
\right)} {|\int_{\gamma(t)} x|^2}
$$
Denote by $W=W_{\gamma,\delta}(\omega, x)$ the determinant in the numerator above.

 It is easily seen that $ W^2$ is
 univalued and hence
rational in $t$. Moreover it has no poles, vanishes at $t=0, -1/4$ and as $t$ tends to infinity it grows no
faster than $t^{(n+1)/2}$. Therefore $W^2$ is a polynomial of degree at most $[(n+1)/2]$ which vanishes at $0$
and $-1/4$, and hence the imaginary part of $F$ along $(-\infty, -1/4)$ has at most $[[(n+1)/2]/2-1]$ zeros.
Summing up the above information we conclude that the increase of the argument of $F(t)$ along the boundary of
$D$ is close to $([(n+1)/2]-1) 2\pi$ or less. Therefore $F$ and hence the Abelian integral $I$ has at most
$[(n+1)/2]-1$ zero in $D$ (and hence in $( -1/4,0)$). It is seen from this proof that the dimension of $V_n$
should be at least $[(n+1)/2]$. Indeed
$$
\int_{\gamma(t)} x^{2k} \equiv 0, \forall k
$$
and
$$
\int_{\gamma(t)} x, \int_{\gamma(t)} x^3, ...
$$
form a basis of $V_n$ (this follows from Proposition \ref{prop3}), which shows that $V_n$ is Chebishev.

2) Suppose now that $\gamma(t) = x_1(t)-x_3(t)$ is a cycle vanishing as $t$ tends to $-1/4$, and defined on the
interval $(-1/4,0)$. The dimension of $V_n$ equals to $n$ and the preceding method does not work. The curve
$\Gamma_f$ is, however, reducible
$$
\Gamma_f =\{ (x,y)\in \C^2: (x+y)(x^2+y^2-1) = 0 \}
$$
and the family $\gamma(t), t\in (-1/4,0)$ corresponds to a piece of the oval $x^2+y^2-1=0$. This oval intersects
$\Gamma_\omega$ in at most $2(n-1)$ points (by Bezout's theorem). The points $(x,y)$ and $(y,x)$ correspond to
$\gamma(t)$  and $-\gamma(t)$ respectively. This shows that each integral $I\in V_n$ can have at most $n-1$
zeros in $(-1/4,0)$. The Proposition is proved.
\end{proof}
\end{exam}
It seems to be difficult to adapt some of the above methods to the case of a general polynomial $f$ of degree
four.

\section{Arithmetic zero dimensional abelian integrals}
\label{s5}
 In this section $\K$ is an arbitrary field of characteristic zero and we work with polynomials in
$\K[x]$. The reader may follow this section for $\K=\Q$. The main result of this section is Theorem \ref{3may05}
which will be used in \S \ref{monodromy} for the functional field $\K=\C(t)$.

For polynomials $f,\omega\in \K[x]$ we define the discriminant of
$f$
$$
\Delta_f:=\prod_{1\leq i,j\leq d}(x_i-x_j)\in\K
$$
and the following polynomial
\begin{equation}
\label{maninnevesht}
\omega*f(x):=(x-\omega(x_1))(x-\omega(x_2))\cdots
(x-\omega(x_d))\in\K[x]
\end{equation}
where
$f(x)=(x-x_1)(x-x_2)\cdots (x-x_d)$.
Note that $(\omega*f)\circ \omega(x_i)=0,\ i=1,2,\cdots,d$ and the
multiplicity of $(\omega*f)\circ \omega$ at $x_i$ is
at least the multiplicity of $f$ at $x_i$. Therefore
\begin{equation}
\label{berlin}
f\mid (\omega*f)\circ \omega
\end{equation}
For $\omega,\omega_1,\omega_2,f, f_1,f_2\in \K[x]$ we have the
following trivial identities:
\begin{equation}
\label{17.6.05}
\omega_1*(\omega_2*f)=(\omega_1\circ \omega_2)*f,\ \omega*(f_1\cdot f_2)=
(\omega*f_1)\cdot (\omega*f_2)
\end{equation}
\begin{prop}
\label{17.6} For an irreducible $f\in\K[x]$ and arbitrary $\omega\in \K[x]$,
we have $\omega*f=g^k$
for some $k\in\N$ and irreducible polynomial $g\in\K[x]$. Moreover,
if for some
simple cycle $\gamma\in H_0(\{f=0\},\Z)$ we have $\int_\gamma \omega=0$
then $k\geq 2$.
\end{prop}
\begin{proof}
Let $$
f=(x-x_1)(x-x_2)\cdots (x-x_d),\ d:=\deg(f), \
I:=\{x_1,x_2,\ldots, x_d\}.
$$
We define
the equivalence relation $\sim$ on $I$:
$$
x_i\sim x_j \Leftrightarrow \omega(x_i)=\omega(x_j)
$$
Let $G_f$ be the Galois group of the splitting field of $f$.
For $\sigma\in G_f$ we have
\begin{equation}
\label{2may2005}
x_i\sim x_j \Rightarrow \sigma(x_i)\sim \sigma(x_j)
\end{equation}
Since $f$ is irreducible over $\K$, the action of $G_f$ on $I$ is
transitive (see for instance \cite{mil05} Prop. 4.4). This  and
(\ref{2may2005}) imply that $G_f$ acts on $I/\sim$ and each
equivalence class of $I/\sim$ has the same number of elements as
others. Let $I/\sim=\{v_1,v_2,\ldots, v_e\},\ e\mid d$ and
$c_i:=\omega(v_i)$. Define
$$
g(x):=(x-c_1)(x-c_2)\cdots (x-c_e)
$$
We have
$$
g^k=f*\omega\in\K[x]
$$
where  $k=\frac{n}{e}$. By calculating the coefficients of $g$ in
terms of the coefficients of the right hand side of the above
equality, one can see easily that $g\in\K[x]$. Since $G_f$ acts
transitively on the roots of $g$, we conclude that $g$ is
irreducible over $\K$.
\end{proof}
Let $f,g,\omega\in \K[x]$ such that
\begin{equation}
\label{2may05}
f\mid g\circ \omega
\end{equation}
We have the morphism
$$
\{f=0\}\stackrel{\alpha_\omega}{\rightarrow}\{g=0\}, \
\alpha_\omega(x)=\omega(x)
$$
defined over $\K$. Let $\gamma\in H_0(\{f=0\},\Z)$ such that
$(\alpha_\omega)_*(\gamma)=0$, where $(\alpha_\omega)_*$ is the
induced map in homology. For instance, if $\deg(g)<\deg(f)$ then
because of (\ref{2may05}), there exist two zeros $x_1,x_2$ of $f$
such that $\int_\gamma \omega=\omega(x_1)-\omega(x_2)=0$ and so the
topological cycle $\gamma:=x_1-x_2$ has the desired property.
Note that and the $0$-form $\omega$ on
$\{f=0\}$ is the pull-back of the $0$-form $x$ by
$\alpha_\omega$.
The
following theorem discusses the inverse of the above situation:
\begin{theo}
\label{3may05} Let $f,\omega\in\K[x]$ be such that such that
\begin{equation}
\label{intvan}
\int_\gamma\omega=0
\end{equation}
for some simple cycle $\gamma\in  H_0(\{f=0\},\Z)$. Then there
exists a polynomial $g\in\K[x]$ such that
\begin{enumerate}
\item
 $\deg(g)< \deg(f)$;
 \item
the degree of each irreducible components of $g$ divides the degree
of an irreducible component of $f$;
\item
$f\mid g\circ \omega$, the morphism
$\alpha_\omega:\{f=0\}\rightarrow \{g=0\}$ defined over $\K$ is surjective and
$(\alpha_\omega)_*(\gamma)=0$.
\end{enumerate}
\end{theo}
The above proposition can be considered as a
zero dimensional counterpart  of the main result of
\cite{evho}  on contraction of curves.
\begin{proof}
Let $f=f_1^{\alpha_1} f_2^{\alpha_2}\cdots f_r^{\alpha_r}$ (resp.
$\omega*f= g_1^{\beta_1} g_2^{\beta_2}\cdots g_s^{\beta_s}$) be the
decomposition of $f\in\K[x]$ (resp. $\omega*f$) into irreducible
components. By Proposition \ref{17.6} and the second equality in
(\ref{17.6.05}), we have $s\leq r$ and we can assume that
$\omega*f_i=g_i^{k_i}$ for $i=1,2,\ldots,s$ and some $k_i\in \N$.
The polynomial $g=g_1^{\alpha_1} g_2^{\alpha_2}\cdots
g_s^{\alpha_s}$ is the desired one.
Except the first item and $(\alpha_\omega)_*(\omega)=0$, all other
parts of the theorem are satisfied by definition.


Let $\gamma=x_1-x_2$. We consider two cases: First let us assume
that $x_1$ and $x_2$ are two distinct roots of an irreducible component of
$f$, say $f_1$. By Proposition \ref{17.6} we have
$\omega*f_1=g_1^{k_1},\ k_1>1$ and so $\deg(g)<\deg(f)$. Now assume
that $x_1$ is a zero of $f_1$ and $x_2$ is a zero of $f_2$. Let
$\omega*f_1=g_1^{k_1},\ \omega*f_2=g_2^{k_2},\ k_1,k_2\in\N$. The
number $\omega(x_1)=\omega(x_2)$ is a root of both $g_i,\ i=1,2$ and
$G_f$ acts transitively on the roots of both $g_i, \ i=1,2$. This
implies that $g_1=g_2$ and so $\deg(g)<\deg(f)$.
\end{proof}

Let $m$ be a prime number and $n<m$.
Theorem \ref{3may05} with $\K=\Q$ implies
that for an irreducible polynomial $f\in \Q[x]$ of degree $m$ the
integral $\int_{\gamma}\omega,\ \omega\in \Q[x]\backslash \Q, \deg(\omega)=n$
never vanishes.
Therefore, the number $Z(m,n,\Q,0)$ cannot be reached by irreducible
polynomials.

\begin{rem}
Let $\K$ be a subfield of $\C$ and $f,\omega\in\K[x]$. Any
$\sigma\in \Gal(\bar\K/ \K)$ induces a map
$$
\sigma: H_0(L_t,\Z)\rightarrow H_0(L_{\sigma(t)},\Z )
$$
in a canonical way and so if $\int_{\gamma(t)}\omega=0$ then
$\int_{\sigma(\gamma(t))}\omega=0$.  This means that
$\Gal(\bar\K/\K)$ acts on the set $\{t\in \C \mid \exists \gamma\in
H_0(L_t,\Z) \hbox{ s.t. } \int_\gamma\omega=0\}$.
\end{rem}
\section{Monodromy group}
\label{s6}
 \label{monodromy} Let $\K$ be a subfield of $\C$, $f\in \K[x]$ and $C$ be the set of its critical
values. We fix a regular value $b$ of $f$.
 The group $\pi_1(\C\backslash C , b)$ acts on $L_b:=\{x_1,x_2,\ldots,x_d\}$
from left. We define the monodromy group
$$
G:=\pi_1(\C\backslash C,b)/\{g\in \pi_1(\C\backslash C,b)\mid
g(x)=x,\ \forall x\in L_b \}\subset S_d
$$
where $S_d$ is the permutation group in $d$ elements
$x_1,x_2,\ldots,x_d$.
 Since the two variable polynomial
$f(x)-t$ is irreducible, the action of $G$ on $L_b$ is also
irreducible. However, the action of $G$ on simple cycles $\Sc\subset
H_0(L_b,\Z)$ may not be irreducible. For instance for $f=x^d,\ b=1$
the group $G$ is generated by the shifting map
$1\mapsto\zeta_d\mapsto\ldots\mapsto\zeta_d^{d-1}\mapsto 1$, where
$\zeta_d=e^{\frac{2\pi i}{2}}$.

Let
$$
\Sc=\Sc_1\cup \Sc_2\ldots\cup \Sc_m
$$
be the partition of $\Sc$ obtained by the action of $G$, i.e. the
partition obtained by the equivalence relation $\gamma_1\sim
\gamma_2$ if $\gamma_1=g\gamma_2$ for some $g\in G$. For
$\omega\in\K[x]$ with $\deg(\omega)=n$ the functions
$$
R_{\omega,i}(t):=\prod_{\gamma\in {\Sc_i/\pm 1}}
\frac{\int_\gamma\omega}{\int_\gamma x},\
\Delta_i(t):=\prod_{\gamma\in\Sc_i} \int_\gamma x,\ i=1,2,\ldots,m
$$
are well-defined in a neighborhood of $b$. They extend to
one valued functions in $\C\backslash C$ and by growth conditions
at infinity and critical values of $f$, we conclude that
they are polynomials in $t$ with coefficients in the
algebraic closure $\bar \K$ of $\K$ in $\C$. Without lose of generality
we assume that for $1\leq m'\leq m$ we have $\Sc_i\not = -\Sc_j$ for
all $1\leq i,j\leq m'$ and for $m\geq i> m'$, $\Sc_i=-\Sc_j$ for some
$1\leq j\leq m'$. Let
$$
R_\omega :=\prod_{i=1}^{m'} R_{\omega,i}=\prod_{\gamma\in {\Sc/\pm 1}}
\frac{\int_\gamma\omega}{\int_\gamma x}
,\ \Delta:=\prod_{i=1}^m\Delta_i=\prod_{\gamma\in\Sc_i} \int_\gamma
x
$$
We define $\tilde f:=f-t$ and consider it as a polynomial in
$\K(t)[x]$. In this way the polynomial $\Delta(t)$ is equal to
$\Delta_{\tilde f}$, the discriminant of $\tilde f$, and
\begin{equation}
\label{18.6.05} R_{\omega}^2=\frac{\Delta_{\omega*\tilde
f}}{\Delta_{\tilde f}}
\end{equation}
 Considering $f=x^d-a_1x^{d-1}-\cdots-a_{d}$ in
parameters $a_i$ with weight$(a_i)=i$, we know that
$\Delta_{\omega*f}$ (resp. $\Delta_f$) is a polynomial (resp.
homogeneous polynomial) of degree $(n-1)d(d-1)$ (resp. $d(d-1)$) in
parameters $a_i$ and with coefficients in $\K$ (resp. in $\Q$). This
implies that we have $R_{\omega},\Delta \in\K[t]$ and by
(\ref{18.6.05})
\begin{equation}
\label{lokht} \sum_{i=1}^{m'}\deg(R_{\omega,i})=\deg(R_{\omega})\leq
\frac{(n-1)(d-1)}{2}
\end{equation}
We summarize the above discussion in the following theorem:
\begin{theo}
\label{main1} Let $\K=\C$ and $\gamma(t)$ be a continuous family of
simple cycles, $\gamma(b)\in \Sc_i$ and  $\omega \in \C[x]$. If the
Abelian integral $ I(t)= \int_{\gamma(t)} \omega $ does not vanish
identically, then the number of its complex zeros in any simply
connected set ${\cal D} \subset \C\setminus \Sigma$  with $b\in
{\cal D}$ is bounded by $\deg(R_{\omega, i})\leq
\frac{(n-1)(d-1)}{2}$.
\end{theo}
The proof follows from the definition of $R_{\omega,i}$ and
(\ref{lokht}). This theorem generalize Theorem \ref{main} and,
roughly speaking, it says that as much as the action of the
monodromy group on $\gamma(b)$ produces less cycles, so far we
expect less zeros for $I(t)$.

In Theorem \ref{main} and \ref{main1} we have assumed that $I(t)$ is
not identically zero. Now the natural question is that if
$I(t)\equiv 0$ then what one can say about $f$ and $\omega$. For
instance, if there is a polynomial $g(x,t)\in \C[t][x]$ such that
$\deg_x(g)<\deg(f)$ and $g(\omega(x),f(x))\equiv 0$ then there is a
continuous family of simple cycles $\gamma(t)$ such that
$\int_{\gamma(t)}\omega\equiv0$.
\begin{theo}
\label{t4}
 Let $\gamma(t)$ be a continuous family of simple cycles in the fibers of $f\in\C[x]$ and $\omega \in
\C[x]$. If the Abelian integral $ I(t)= \int_{\gamma(t)} \omega $ vanishes identically, then there is a
polynomial $g(x,t)\in\C[t][x]$ such that
\begin{enumerate}
\item
$\deg_x(g)\mid \deg(f)$ and $\deg_x(g)<\deg(f)$;
\item
$g(\omega(x),f(x))\equiv 0$;
\item
If the action of the monodromy group on a regular fiber of $f$ is
irreducible then $\omega=p(f)$ for some $p\in\C[x]$.
\end{enumerate}
\end{theo}
\begin{proof}
We consider $\tilde f=f-t$ as a polynomial in $\K[x]$ with
$\K=\C(t)$. The assumption of the theorem is translated into
$\int_{\gamma}\omega=0$ for some simple cycle $\gamma\in
H_0(\{\tilde f=0\},\Z)$. We apply Theorem \ref{3may05} and we
conclude that there are polynomials $g,s\in \K[x],\
\deg(g)<\deg(f),\ \deg(g)\mid \deg(f)$ such that $s\cdot
(f-t)=g(\omega(x),t)$. Note that $\tilde f$ is irreducible over
$\K$. After multiplication with a certain element in $\C[t]$, we can
assume that $s,g\in \C[x,t]$. We replace $t$ with $f(x)$ and in this
way the item 1 and 2 are proved.

The third part of the theorem follows from Proposition \ref{prop2}.
We give an alternative proof as follows: We identify the elements of
the splitting field of $\tilde f$ with holomorphic functions
in a neighborhood of $b$ in $\C$. In this way we can identify the
monodromy group $G$ with a subset of the Galois group $G_{\tilde f}$
of the splitting field of $\tilde f$ over $\K$. If the action of $G$
on simple cycles is irreducible then by Theorem \ref{3may05} (in fact
its proof) we have $g=x-p(t)$ for some $p(t)\in\C[t]$.
\end{proof}

The space of polynomials with $d-1$ distinct
critical values can be identified with a quasi-affine subset $T$ of $\C^d$.
We claim that for $f\in T$, the action of the
monodromy group $G$ on $\Sc$ is irreducible. Since our assertion is
topological and $T$ is connected, it is enough to prove our
assertion for an example of $f\in T$; for instance take an small
perturbation $\tilde f$ of $f=(x-1)(x-2)\cdots(x-d)$ which has $d-1$
non-zero distinct critical values $\tilde c_1,\tilde
c_2,\ldots,\tilde c_{d-1}$. Let $b=0$. We take a system of distinguished
paths $s_i,\ i=1,2,\ldots,d-1$ in $\C$ (see
\cite{arn}) such that $\gamma_i=\tilde i-\tilde{(i+1)},\ 1\leq i\leq d-1$
vanishes along $s_i$ in the critical point associated to
$\tilde c_i$, where $\tilde i\in \tilde f^{-1}(0)$ is near $i\in
f^{-1}(0)$. Now, the intersection graph of $\gamma_i$'s (known as
Dynkin diagram of $f$) is a line graph, and so it is connected. By
Picard-Lefschetz formula in dimension zero we conclude that the
action of the monodromy group on simple cycles of $H_0(\tilde
f^{-1}(b),\Z)$ is irreducible.



\end{document}